\documentclass{amsart}
\usepackage{amssymb}
\usepackage{amsmath}
\usepackage{amsthm}
\usepackage[ansinew]{inputenc}
\usepackage[all]{xy}
\xyoption{2cell}

\newcommand{\F}{{\mathcal F}}
\newcommand{\Ol}{{\mathcal O}}
\newcommand{\f}{\varphi}

\newcommand{\ac}{{\mathcal H}}
\newcommand{\E}{{\mathcal E}}
\newcommand{\hD}{{\widehat D}}
\newcommand{\pu}{{\mathbb P^1}}
\newcommand{\proj}{\mathbb P}

\newcommand{\pd}{{\mathbb P^2}}

\newcommand{\tl}{\widetilde}
\newcommand{\Dcap}{\widehat D}

\newcommand{\G}{{\mathbb G}}
\newcommand{\Z}{{\mathbb Z}}
\newcommand{\qu}{\mathcal Q}

\newcommand{\lra}{{\longrightarrow}}
\newcommand{\shse}[3]{0 ~\lra ~#1~ \lra ~#2~ \lra ~#3~ \lra~ 0}
\DeclareMathOperator{\loc}{\mathrm{Locus}}

\newcommand{\ratcurves}{\textrm{Ratcurves}^n(X)}

\DeclareMathOperator{\cone}{NE}
\DeclareMathOperator{\pic}{Pic}

\DeclareMathOperator{\Exc}{Exc}

\newcommand{\twospan}[2]{\langle #1,#2 \rangle}

\newcommand{\ba}{{B\v anic\v a~}}
\newcommand{\rc}[2]{#1 \xymatrix{\ar@{-->}[r] & }{#2}}

\newtheorem{theorem}{Theorem}[section]
\newtheorem{lemma}[theorem]{Lemma}
\newtheorem{proposition}[theorem]{Proposition}
\newtheorem{corollary}[theorem]{Corollary}

\newtheorem{teo}{Theorem}[subsection]
\newtheorem{lem}[teo]{Lemma}
\newtheorem{prop}[teo]{Proposition}
\newtheorem{cor}[teo]{Corollary}

\theoremstyle{definition}

\newtheorem{example}[theorem]{Example}
\theoremstyle{remark}
\newtheorem{remark}[theorem]{Remark}
\newtheorem{question}[theorem]{Question}

\theoremstyle{definition}
\newtheorem{defi}[teo]{Definition}
\newtheorem{ex}[teo]{Example}
\theoremstyle{remark}

\numberwithin{equation}{section}

\begin{document}

\author{Carla Novelli}
\address{Dipartimento di Matematica,\newline Universit\`a di Genova,\newline via Dodecaneso 35,\newline I-16146 Genova}
\curraddr{Dipartimento di Matematica e Applicazioni, \newline Universit\`a di Milano - Bicocca,\newline via R. Cozzi 53, \newline I-20126 Milano}
\email{carla.novelli@unimib.it}
\author{Gianluca Occhetta}
\address{Dipartimento di Matematica,\newline Universit\`a di Trento, \newline via Sommarive 14,\newline I-38050 Povo (TN)}
\email{gianluca.occhetta@unitn.it}

%\keywords{Fano manifolds, vector bundles, extremal rays, rational curves}

\title[Manifolds containing a linear subspace]{Projective manifolds containing a large linear subspace with nef normal bundle}

%\begin{abstract}
%We classify smooth complex projective varieties $X \subset \proj^N$ of dimension $2s+1$
%containing a linear subspace $\Lambda$ of dimension $s$ whose normal bundle
%$N_{\Lambda/X}$ is numerically effective.
%\end{abstract}

\maketitle

\section{Introduction}

Let $X \subset \proj^N$ be a smooth complex projective variety of dimension $n$,
containing a linear subspace $\Lambda$ of dimension $s$; denote by $N_{\Lambda/X}$
its normal bundle and by $c$ the degree of the first Chern class of $N_{\Lambda/X}$. \\
If $N_{\Lambda/X}$ is numerically effective, then $X$ is covered by lines; if
furthermore $s$ is sufficiently large, the restrictions imposed on $X$ become
stronger.\par
\medskip
By \cite[Theorem 2.5]{BSW}, if $s+c > \frac{n}{2}$ then, for large $m$, denoting by $H$ the restriction to $X$
of the hyperplane bundle, the linear system
$|m(K_X + (s+1+c)H)|$ defines an extremal ray contraction of $X$ which contracts $\Lambda$.
If $s$ is greater than $\frac{n}{2}$ this contraction is a projective bundle, as shown in \cite{SatoTo}
(see also \cite[Theorem 2.5]{BI}).
The same result holds if $s=\frac{n}{2}$ and  $N_{\Lambda/X}$ is trivial,
by  \cite[Theorem 1.7]{Ein} and \cite[Theorem 2.4]{Widef}.\par
\medskip
The complete study of the case $s=\frac{n}{2}$ is the subject of \cite{SatoTo};
the setup of the quoted paper is different - it is not assumed the existence of a linear
space of dimension $\frac{n}{2}$ with nef normal bundle, but the existence of a linear
space of dimension $\frac{n}{2}$ through every point of $X$ -  yet the assumptions
are in fact equivalent.\\
The most difficult cases in \cite{SatoTo} are manifolds of Picard number one, which turn
out to be, besides linear spaces, hyperquadrics and Grassmannians of lines.\par
\medskip
In this paper we study the next case, {\em i.e.} $n=2s+1$, proving the following

\begin{theorem}\label{main}  Let $X \subset \proj^N$ be a smooth variety of dimension $2s+1$,
containing a linear subspace $\Lambda$ of dimension $s$, such that its normal bundle $N_{\Lambda/X}$ is numerically
effective. If the Picard number of $X$ is one, then $X$ is one of the following:
\begin{enumerate}
\item a linear space $\proj^{2s+1}$;
\item a smooth hyperquadric $\mathbb Q^{2s+1}$;
\item a cubic threefold in $\proj^4$;
\item a complete intersection of two hyperquadrics in $\proj^5$;
\item the intersection of the Grassmannian of lines $\G(1,4) \subset \proj^9$ with three general hyperplanes;
\item a hyperplane section of the Grassmannian of lines $\G(1,s+2)$ in its Pl\"ucker embedding.
\end{enumerate}
If the Picard number of $X$ is greater than one, then there is an elementary
contraction $\f \colon X \to Y$ which contracts $\Lambda$ and one of the following occurs:
\begin{enumerate}
\item[(7)] $\f \colon X \to Y$ is a scroll;
\item[(8)] $Y$ is a smooth curve, and the general fiber of $\f$ is
\begin{itemize}
\item[(8a)] the Grassmannian of lines $\G(1,s+1)$;
\item[(8b)] a smooth hyperquadric $\mathbb Q^{2s}$;
\item[(8c)] a product of projective spaces $\proj^s \times \proj^s$.
\end{itemize}
\end{enumerate}
\end{theorem}

The outline of the paper is the following: first of all, we use the theory of uniform vector bundles
on the projective space, together with some standard exact sequences, to classify all possible normal bundles
$N_{\Lambda/X}$.
Then we consider separately the case of Picard number greater than one and the case of Picard number one;
in fact the ideas and the proofs are very different.\par
\medskip

If the Picard number is greater than one, we combine the ideas and techniques of \cite{BSW}
with those of \cite{BCD} to show that a dominating family of lines on $X$ of anticanonical degree
$\ge \frac{n+1}{2}$ is extremal, {\em i.e.} the numerical class of a line spans a Mori extremal ray  of
$\cone(X)$.
The contraction of this ray is the morphism $\f\colon X \to Y$ appearing in the second part of the
statement of Theorem (\ref{main}). The general fiber $F$ of $\f$ is then a manifold covered by
linear spaces of dimension $\ge \frac{\dim F}{2}$, and this leads to its classification.\par
\medskip
If the Picard number is one, the main idea is to study the manifold $\tl X$ obtained by blowing-up $X$
along $\Lambda$; we prove that $\tl X$ is a Fano manifold, and then we study its ``other" extremal contraction.
As a first application of this construction, in section (5) we show how to use it to complete
\cite[Main Theorem]{SatoTo}.\\
In the setup of Theorem (\ref{main}), the hardest case corresponds
to the normal bundle $N_{\Lambda/X} \simeq T_\Lambda(-1) \oplus \Ol_\Lambda$,
which gives rise to case $(6)$. In this case we need to use twice the blow-up
construction: first we blow-up $X$ along $\Lambda$ and we show that there is a special
one-parameter family $\Sigma$ of linear spaces to which $\Lambda$ belongs; then we blow-up
$X$ along $\Sigma$ and, studying this blow-up, we are able to describe completely the variety.

\section{Background material}

%\begin{definition}
A smooth complex projective variety $X$ is called {\em Fano}
if its anticanonical bundle $-K_X$ is ample; the {\em index} $r_X$ of $X$
is the largest natural number such that $-K_X=mH$ for some (ample) divisor
$H$ on $X$. Since $X$ is smooth, $\pic (X)$ is torsion free, therefore the
divisor $L$ satisfying $-K_X=r_XL$ is uniquely determined and called the {\em fundamental divisor} of $X$.
Fano manifolds with $r_X=\dim X -1$ are called {\em del Pezzo} manifolds.
%\end{definition}

\subsection{Extremal contractions}

Let $X$ be a smooth projective variety of dimension $n$ defined over the field of complex numbers.\\
A {\em contraction} $\f \colon X \to Z$ is a proper surjective map with connected fibers onto a normal variety $Z$.\\
If the canonical bundle $K_X$ is not nef, then the negative part of
the cone $\cone(X)$ of effective 1-cycles is locally polyhedral, by the Cone Theorem.
By the Contraction Theorem, to every face in this part of the cone is associated
a contraction, called {\em extremal contraction}
or {\em Fano--Mori contraction}.\\
An extremal contraction associated to a face of dimension one, {\em i.e.} to an extremal ray, is called an
{\em elementary contraction}.\\
A Cartier divisor $H$ such that $H = \f^\ast A$ for an ample divisor $A$ on $Z$
is called a {\em supporting divisor} of the contraction $\f$.

\begin{defi} An elementary fiber type extremal contraction
$\f\colon X \to Z$ is called a {\em scroll} (respectively a {\em quadric fibration})
if there exists a $\f$-ample line bundle $L \in \pic(X)$ such that
$K_X+(\dim X-\dim Z+1)L$ (respectively $K_X+(\dim X-\dim Z)L$)
is a supporting divisor of $\f$.\\
An elementary fiber type extremal contraction $\f\colon X \to Z$ onto a smooth
variety $Z$ is called a $\mathbb P${\em-bundle} (respectively
{\em quadric bundle}) if there exists a vector bundle $\E$ of rank
$\dim X-\dim Z+1$ (respectively of rank $\dim X-\dim Z+2$) on $Z$ such that
$X \simeq \proj_Z(\E)$ (respectively there exists an embedding of $X$ over $Z$ as a divisor of
$\proj_Z(\E)$ of relative degree 2).\\
%An equidimensional scroll is a $\proj$-bundle by \cite[Lemma 2.12]{Fuj4},
%while an equidimensional quadric fibration is a quadric bundle by
%\cite[Theorem B]{ABWtwo}.\\
Some special scroll contractions arise from projectivization of
B\v anic\v a sheaves (cf. \cite{BalW}); in parti\-cu\-lar,
if $\f\colon X \to Z$ is a scroll such that every fiber has dimension $\le \dim X- \dim Z+1$,
then $Z$ is smooth and $X$ is the projectivization of a B\v anic\v a sheaf on $Z$
(cf. \cite[Proposition 2.5]{BalW}); we will call these contractions
{\em special B\v anic\v a scrolls}.
\end{defi}

\subsection{Families of rational curves}

Let $X$ be a smooth projective variety of dimension $n$ defined over the field of complex numbers.

\begin{defi} \label{Rf}
A {\em family of rational curves} is an irreducible component
$V \subset \ratcurves$ (see \cite[Definition 2.11]{Kob}).
Given a rational curve, we will call a {\em family of deformations} of that curve
any irreducible component of $\ratcurves$ containing the point parametrizing that curve.\\
We will say that $V$ is {\em unsplit} if it is proper.\\
We define $\loc(V)$ to be the set of points of $X$ through which there is a curve among those
parametrized by $V$ and we say that $V$ is a {\em dominating family} if $\overline{\loc(V)}=X$.\\
We denote by $V_x$ the subscheme of $V$ parametrizing rational curves passing through $x \in \loc(V)$
and by $\loc(V_x)$ the set of points of $X$ through which there is a curve among those
parametrized by $V_x$.\\
By abuse of notation, given a line bundle $L \in \pic(X)$, we will denote by $L \cdot V$
the intersection number $L \cdot C_V$, with $C_V$ any curve among those
parametrized by $V$.
\end{defi}

\begin{defi}
An unsplit dominating family $V$ defines a relation of rational connectedness with respect
to $V$, which we shall call {\em rc$(V)$-relation} for short, in the following
way: $x$ and $y$ are in rc$(V)$-relation if there
exists a chain of rational curves among those
parametrized by $V$ which joins $x$ and $y$.
\end{defi}

To the rc$(V)$-relation we can associate a fibration, at least on an open subset
(\cite{Cam81}, \cite[IV.4.16]{Kob}); we will call it {\em rc$(V)$-fibration}.

\begin{prop}\cite[IV.2.6]{Kob}\label{iowifam}
Let $V$ be an unsplit family of rational curves on $X$. Then
 \begin{itemize}
      \item[(a)] $\dim X -K_X \cdot V \le \dim \loc(V)+\dim \loc(V_x) +1$;
      \item[(b)] $-K_X \cdot V \le \dim \loc(V_x)+1$.
   \end{itemize}
\end{prop}

This last proposition, in case $V$ is the unsplit family of deformations of a minimal extremal
rational curve, {\em i.e.} of a rational curve of minimal anticanonical degree in an extremal face of
$\cone(X)$, gives the {\em fiber locus inequality}:

\begin{prop}\cite[Theorem 0.4]{Io}, \cite[Theorem 1.1]{Wicon}\label{fiberlocus}
Let $\f$ be a Fano--Mori contraction of $X$.
Denote by $E$ the exceptional locus of $\f$ and
by $F$ an irreducible component of a non-trivial fiber of $\f$. Then
$$\dim E + \dim F \geq \dim X + l -1,$$
where $l :=  \min \{ -K_X \cdot C\ |\  C \textrm{~is a rational curve in~} F\}$.
If $\f$ is the contraction of an extremal ray $R$, then $l(R):=l$ is called the {\em length of the ray}.
\end{prop}

\begin{defi}
Let $V$ be an unsplit family of rational curves
on $X$ and $Z \subset X$. We denote by $\loc(V)_Z$ the set of points $x \in X$ such that there exists
a curve $C$ in $V$ with $C \cap Z \not = \emptyset$ and $x \in C$.
\end{defi}

We will use some properties of $\loc(V)_Z$, summarized in the following

\begin{lem}\label{locvf}\cite[Section 2]{CO}, \cite[Proof of Lemma 1.4.5]{BSW}
Let $Z \subset X$ be a closed subset and $V$ an unsplit family.
Assume that curves contained in $Z$ are numerically independent from curves in $V$, and that
$Z \cap \loc(V) \not= \emptyset$. Then
$$\dim \loc(V)_Z \ge \dim Z -K_X \cdot V - 1.$$
If $\sigma$ is an extremal face of $\cone(X)$, $F$ is a fiber
of the contraction associated to $\sigma$ and $V$ is an unsplit family, numerically independent from curves
whose numerical class is in $\sigma$, then
$$\cone(\loc(V)_F, X) = \langle \sigma, [V] \rangle,$$
{\em i.e.} the numerical class in $X$ of a curve in $\loc(V)_F$ is in the subcone of $\cone(X)$
generated by $\sigma$ and $[V]$.
\end{lem}

\subsection{Some extremal contractions related to Grassmannians}

We will now present some examples of Fano manifolds admitting a projective bundle structure
and another extremal contraction $\f$ whose target is a Grassmannian of lines.
We will use these descriptions later in our proofs.
%We will always denote by $\G(1,s)$ the Grassmannian of lines
%in a projective space of dimension $s$.

\begin{ex}\label{ex1} Let $\G(1,s)$ be the Grassmannian of lines in $\proj^s$ and denote by
$\mathcal I$ the incidence variety. Consider the incidence diagram:
$$\xymatrix{&\mathcal I \ar[rd]^\f \ar[ld]_p& \\ \proj^s & & **{!<0.2cm,0cm>}{\G(1,s)}.}$$
Then $p$ and $\f$ are projective bundles, namely $\mathcal I =\proj_{\proj^s}(p_\ast \f^\ast \Ol_{\G(1,s)}(1))=
\proj_{\proj^s}(\Omega_{\proj^s}(2))$ and $\mathcal I =\proj_{\G(1,s)}(\f_\ast p^\ast \Ol_{\proj^s}(1))=
 \proj_{\G(1,s)}(\qu)$, where $\qu$ is the universal quotient bundle on $\G(1,s)$.
\end{ex}
%<1ex>
\begin{ex}\label{qb}
As in the previous example, let $\G(1,s)$ be the Grassmannian of lines in $\proj^s$ and let
$\mathcal I$ denote the incidence variety; consider the following diagram, obtained by the incidence diagram:
$$\xymatrix@C=10pt{  & \mathcal  I \times \pu  \ar[rd]^q \ar[ld]_p \ar@/^1.2pc/[rrd]^(.55)\f & \\
  **{!<-0.2cm,0cm>}{\proj^s \times \pu}  & &  **{!<0.4cm,0cm>}{\G(1,s)\times \pu}  \ar[r]^(.50)g & \G(1,s) .}$$
%$$\xymatrix@C=25pt{ \mathcal  I \times \pu  \ar[r]^q \ar[d]_p & \G(1,s)\times \pu \ar[r]^g \ar[d]
%& \G(1,s) \\ \proj^s \times \pu \ar[d]_{p_1} \ar[r]^{p_2} & \pu & \\
%\proj^s & & }$$
The composition $\f= g \circ q$ gives a morphism $\f \colon \mathcal I \times \pu \to \G(1,s)$ whose
fibers are smooth two-dimensional quadrics.
Let $\ac$ be $p^\ast\Ol_{\proj^s \times \pu}(1,1)$ and put $\E:=\f_\ast \ac$.
We have
$$\E=\f_\ast \ac= g_\ast(q_\ast \ac)=g_\ast(\Ol_{\G(1,s)\times \pu}(1,0) \otimes g^\ast \qu)=\qu^{\oplus 2}.$$
The product $\mathcal I \times \pu = \proj_{\proj^s \times \pu}(p_1^\ast \Omega_{\proj^s}(2))$,
where $p_1$ denotes the projection onto $\proj^s$,
embeds in $\proj_{\G(1,s)}(\E)$ as a divisor of relative degree $2$, {\em i.e.} it belongs to a linear system
$| 2\ac - \f^\ast L|$ for some line bundle $L$ in $\pic(\G(1,s))$.
The discriminant divisor of the quadric bundle is in the linear system $|2\det \E -4L|$
and it is trivial, since every fiber of $\f$ is smooth. It follows that $L =\Ol_{\G(1,s)}(1)$.
\end{ex}

\begin{ex}\label{gblow} Let $\G(1,s+1)$ be the Grassmannian of lines in $\proj^{s+1}$ and let
$\G(1,H) \subset \G(1,s+1)$ be a sub-Grassmannian corresponding to the lines of $\proj^{s+1}$
contained in a fixed hyperplane $H$.\\
Consider the rational map $\psi \colon \rc{\G(1,s+1)}{H}$, which associates to a line $l$
the point of intersection of $l$ with $H$. This map is not defined precisely along
the points of $\G(1,s+1)$ representing the lines contained in $H$,
{\em i.e.} along the sub-Grassmannian $\G(1,H)$.\\
Consider the resolution of $\psi$, obtained by blowing-up $\G(1,s+1)$ along $\G(1,H)$:
$$\xymatrix@C=45pt@R=35pt{\tl \G(1,s+1)   \ar[d]_{p}  \ar[r]^(.43)\f & **{!<0.4cm,0cm>}{\G(1,s+1)} \ar@{-->}[ld]^\psi \\
H & }$$
The contraction $p$ is a $\proj^s$-bundle over $H$, whose fibers
are the strict transforms of linear subspaces $\proj^s \subset \G(1,s+1)$ corresponding
to stars of lines with center in $H$, namely
$\tl \G(1,s+1)=\proj_H(p_\ast \f^\ast \Ol_{\G(1,s+1)}(1))= \proj_H(\Omega_{\proj^s}(2) \oplus \Ol_{\proj^s}(1))$.
\end{ex}

\section{Manifolds with Picard number greater than one}

In this section we are going to show that, if the Picard number of $X$ is greater than one, and $X$ is covered by linear spaces of dimension $s \ge [n/2]$, then, either $\pic(X) \simeq \mathbb Z$ or there is an elementary Mori contraction to a positive dimensional variety, whose general fiber is covered by linear spaces; in this last case we will then get the description of the general fiber by Corollary (\ref{exsato}).\par 
\medskip
Let $X$ be a smooth complex projective variety, let $V$ be an unsplit dominating family of rational curves for $X$
and let $q\colon \rc{X}{Y}$ be the rc$(V)$-fibration.
Let $B$ be the indeterminacy locus of $q$; notice that $\dim B \le \dim X-2$, as $X$ is smooth.
Moreover, by \cite[Proposition 1]{BCD} $B$ is the union of all rc$(V)$-equivalence classes
of dimension greater than $\dim X - \dim Y$.

\begin{lemma}\label{inB}
Let $V$ be an unsplit dominating family of rational curves on a smooth projective variety $X$.
Let $B$ be the indeterminacy locus of the rc$(V)$-fibration $q\colon \rc{X}{Y}$,
let $D$ be  very ample  on $q(X \setminus B)$ and let $\widehat D:=\overline{q^{-1}D}$. Then
\begin{itemize}
\item[(a)] $\widehat D \cdot V=0$;
\item[(b)] if $C \not \subset B$ is a curve whose numerical class is not proportional to $[V]$, then $\Dcap \cdot C >0$;
\item[(c)] if $[V]$ does not span an extremal ray of $\cone(X)$, then there exists
a curve $C \subset B$ whose class is not proportional to $[V]$ such that $\widehat D \cdot C \le 0$.
\end{itemize}
\end{lemma}

\begin{proof}
A general cycle of $V$ is contained in a fiber of $q$ disjoint from $\widehat D$,
so $\widehat D \cdot V=0$.\\
If $C$ is as in (b), then $q(C)$ is a curve in $Y$ and the result follows from
projection formula.\\
Finally, if $[V]$ does not span an extremal ray, then either $\hD$ is not nef or
$\hD$ is nef but $\hD^{\bot} \cap \cone(X) \supsetneq [V]$.
In both cases there exists a curve $C \subset X$ whose class is not
proportional to $[V]$ such that $\widehat D \cdot C \le 0$. Such a curve must be contained in $B$
by \cite[Proof of Proposition 1]{BCD}.
\end{proof}

\begin{lemma}\label{pippo} Let $X$ be a manifold which admits an unsplit dominating family of
rational curves~$V$. Assume that there exists an extremal face
$\Sigma \subseteq \overline{\cone(X)}_{K_X <0}$ such that $[V] \subset \Sigma$.\\
Then, either $[V]$ spans an extremal ray, or there exists an extremal ray in $\Sigma$
whose exceptional locus is contained in the indeterminacy locus $B$ of the rc$(V)$-fibration.
In particular, this ray is associated with a small contraction.
\end{lemma}

\begin{proof}
Let $\tau$ be a minimal subface of $\Sigma$ containing $[V]$. If $\dim \tau=1$, then $[V]$ spans an extremal ray.\\
Assume that $\dim \tau \ge 2$. Let $\widehat D$ be as in Lemma (\ref{inB}). Since $\widehat D \cdot V=0$, then
either $\widehat D$ is zero on every extremal ray of $\sigma$ or it is negative on at least one ray.
In both cases, by part (b) of Lemma (\ref{inB}) there is at least one ray whose exceptional
locus is contained in $B$, and the assertion follows as $\dim B \le \dim X-2$.
\end{proof}

The following is a slight improvement of \cite[Theorem 2.5]{BSW} (cf. also \cite[Theorem 2.4]{BI},
where the case $-K_X \cdot V  \ge \frac{n + 3}{2}$ is treated):

\begin{theorem}\label{extremal}
Let $(X,H)$ be a polarized $n$-fold with a dominating family of rational curves $V$ such that $H \cdot V=1$.
If $-K_X \cdot V  \ge \frac{n + 1}{2}$, then~$[V]$ spans an extremal ray of $\cone(X)$.
\end{theorem}

\begin{proof} Denote by $m$ the positive integer $-K_X \cdot V$ and by $L$ the adjoint divisor $K_X +mH$.\par
\medskip
Assume first that $L$ is nef.\\
Denote by $q \colon \rc{X}{Y}$ the rc$(V)$-fibration and by $B$ its indeterminacy locus.
Assume that $[V]$ does not span an extremal ray in $\cone(X)$.
This implies that $L$ defines an extremal face $\Sigma$ of dimension at least two, containing $[V]$.\\
By Lemma (\ref{pippo}) there exists an extremal ray $R \in \Sigma$
whose associated contraction $\f$ is small; moreover, since $L \cdot R = 0$
the length of this extremal ray is greater than or equal to $m$.\\
If $F$ is a non-trivial fiber of $\f$, by Proposition (\ref{fiberlocus}), we have $\dim F \ge m+1$.\\
Let $x$ be a point in $F$; $\loc(V_x)$ meets $F$, but, since $[V]$ is independent
from $R$, the intersection has to be zero-dimensional. This implies that
$$\dim \loc(V_x) \le n-m-1 \le \frac{n-3}{2},$$
contradicting part (b) of Proposition (\ref{iowifam}).\par
\medskip
Assume now that $L$ is not nef.\\
This assumption yields the existence of an extremal ray
$R$ such that $L \cdot R <0$.
Notice that $R$ has length $\ge m+1$, hence every non-trivial fiber of the associated contraction
has dimension $\ge m$ by Proposition (\ref{fiberlocus}).\\
We have, by Lemma (\ref{locvf}),
$$\dim X \ge \dim \loc(V)_F \ge -K_X \cdot V + \dim F -1\ge m +m -1 \ge n,$$
hence $\loc(V)_F =X$. We can apply the second part of Lemma (\ref{locvf})
to get $\cone(X)= \twospan{[V]}{R}$ and we are done.
\end{proof}

\begin{remark} Combining the ideas and tecniques of \cite{BSW} with those of \cite{BCD} it is
actually possible to prove the statement of Theorem (\ref{extremal}) under the weaker assumption
that $-K_X \cdot V  \ge \frac{n - 1}{2}$; however the proof becomes very long and complicated, so,
since it is not necessary for our main theorem, we will
present it elsewhere (\cite{NOlines}).
\end{remark}

\begin{theorem}\label{nc} Let $X \subset \proj^N$ be a smooth variety of dimension $n$, covered by linear spaces of dimension $s \ge [n/2]$; then there is an elementary Mori contraction $\f \colon X \to Y$, which contracts lines in the corresponding covering family. Moreover, either $\pic(X) \simeq \mathbb Z$ and $Y$ is a point, or, denoting by $F$ a general fiber of $\f$, one of the following occurs: $\pic(F) \simeq \mathbb Z$ or $n=2s+1$ and $F \simeq \proj^s \times \proj^s$.
\end{theorem}

\begin{proof} Let  $l$ be a general line in a general linear space;  by the assumptions there is a dominating family of lines in $X$ containing $l$. By adjunction we have $-K_X \cdot l \ge s+1$, hence, by Theorem
(\ref{extremal}), the numerical class of $l$ spans an extremal ray of $X$.\\
Let $\f\colon X \to Y$ be the contraction of this extremal ray and
let $F$ be a general fiber of $\f$; $F$ has dimension at most $2s$ and, by adjunction, is a
Fano manifold of index $\ge s+1$, hence either $F \simeq \proj^s \times \proj^s$ or $\rho_F=1$ by
\cite[Theorem B]{Wimu}.
\end{proof}

\begin{example} We show with an example that the last case of Theorem (\ref{nc}) is effective;
the idea on which it is based has been suggested by Wi\'sniewski for \cite[Example 7.2]{Opa}.\\
Let $C'$ be a smooth curve with a free $\Z_2$-action, so that the action induces an \'etale
covering $\pi \colon C' \to C$ of degree $2$. Let $G$ be $\proj^s \times \proj^s$
and take on $G$ the $\Z_2$-action which exchanges the factors.\\
Let $X':=G \times C'$ and denote by $X$ the quotient of $X'$ by the product action of $\Z_2$;
the action is free and so $X$ is smooth.
By the universal property of group actions there exists a morphism $\f\colon  X \to C$
such that the following diagram commutes:
$$
\xymatrix{
X' \ar[rr] \ar[d]_{\pi'} &&  C' \ar[d]^{\pi}\\
X \ar[rr]^{\f} && C
}
$$
The map $\f\colon X \to C$ is an extremal contraction and every fiber is a product of projective spaces
$\proj^s \times \proj^s$. We will now show that $\f$ is elementary.\\
Let $l$ be a line in  $G$ and consider the product
$l \times C' \subset G \times C' =X'$: it is a flat family of rational curves.
Let $c$ be a point of $C$ and let  $\{c_1',c_2'\}$ be $\pi^{-1}(c)$; finally set
$l_i' := l \times \{c_i'\}$ and consider the restriction of the previous diagram

$$\xymatrix{
G \times \{c_1',c_2'\} \ar[rr] \ar[d]_{\pi'} &&   \{c_1',c_2'\}\ar[d]^{\pi}\\
\f^{-1}(c) \simeq \proj^s \times \proj^s \ar[rr]^(.60){\f} && c
}$$
Since the product action identifies $G \times \{c_1'\}$
with $G \times\{c_2'\}$ exchanging the factors
we have that $l_1= \pi'(l_1')$ is a line in a fiber of
the projection of $\f^{-1}(c)$ onto the first factor and $l_2=\pi'(l_2')$
is a line in a fiber of the projection of $\f^{-1}(c)$ onto the second factor, hence
lines in the two factors are algebraically and thus numerically equivalent.
\end{example}

\begin{remark} In the last example, $X$ has an unsplit dominating family of rational curves $V$
such that $V_x$ has dimension $\frac{\dim X-3}{2}$ and is reducible for every $x$.\\
This should be compared with  \cite[Theorem 5.1]{KeKo}, in which it is proved that, if
$\dim V_x \ge \frac{\dim X-1}{2}$, then $V_x$ is irreducible.
\end{remark}

\section{A general construction}

In this section we will present a blow-up construction and we will show how to apply it to manifolds of
Picard number one containing a linear space with nef normal bundle.
\\
The construction in the following proposition has been inspired by the graduate thesis
\cite{TesiPa} supervised by the second named author.

\setcounter{equation}{0}

\begin{proposition}\label{gc}
Let $X \subset \proj^N$ be a Fano manifold  of dimension $n$, index $r_X$
and Picard number one, covered by lines and containing a smooth subvariety $\Sigma$ of dimension $s$
which is the intersection of its linear span with $X$,
{\em i.e.} $\Sigma = X \cap \langle \Sigma \rangle$,
with $[n/2] \le s \le n-2$.
Let $\sigma \colon \tl X \to X$ be the blow-up of $X$ along $\Sigma$ and let
$E=\proj_{\Sigma}(N^\ast_{\Sigma/X})$ be the exceptional divisor of $\sigma$. \\
Then $\cone(\tl X)=\twospan{[C_\sigma]}{[{\ell}]}$, where $C_\sigma$ is a minimal
curve contracted by $\sigma$ and $\ell$ is the strict transform of a line meeting $\Sigma$ at one point. \\
If $r_X \ge [n/2]+1$, then $\tl X$ is a Fano manifold, the length of the ray $\mathbb R_+[\ell]$ is $r_X-n+s+1$ and the extremal contraction $\f \colon \tl X \to Y$ associated to $\mathbb R_+[\ell]$ is  the morphism  defined by the linear system $|m(\sigma^\ast\Ol_X(1)-E)|$ for $m>>0$.\\
If $r_X \ge [(n+1)/2]+1$, then also $E$ is a Fano manifold;  moreover, for any positive $m$, the restriction to $E$ of  the morphism given by $|m(\sigma^\ast\Ol_X(1)-E)|$ is the morphism given by the linear system $|m\xi_{N^\ast_{\Sigma/X}(1)}|$.\end{proposition}

\begin{proof}
Since the Picard number of $X$ is one and $X$ is covered by lines, it follows by
\cite[Proposition 1.1]{AWinv} that $X$ is rationally connected with respect to a dominating family of lines.\\
Consider the rational map $\rc{X}{\tl Y}$ defined by the linear system
$|\Ol_X(1) \otimes {\mathcal I}_\Sigma|$ of hyperplanes containing $\Sigma$.
Let $\tl \f \colon \tl X \to \tl Y$ be the resolution of this map.
Then the morphism $\tl \f$ is defined by the linear system $|\ac -E|$, where $\ac$ denotes  the pull-back
$\sigma^\ast\Ol_X(1)$.\\
Let $l \subset X$ be a line meeting $\Sigma$ but not contained in it;
notice that such a line exists because $X$ is rationally connected with respect to a family of lines.\\
Since $\Sigma = X \cap \langle \Sigma \rangle$, the intersection $l \cap \Sigma$ consists of one point,
hence the morphism $\tl \f$ contracts $\ell$, the strict transform of $l$. Therefore,
denoting by $C_\sigma$ a rational curve of minimal degree
contracted by $\sigma$, we obtain $\cone(\tl X)=\twospan{[C_\sigma]}{[{\ell}]}$.
The contraction associated to the ray $R=\mathbb R_+[\ell]$ is therefore given by the Stein
factorization of $\tl \f$, {\em i.e.} it is defined by the linear system $|m(\ac -E)|$ for $m >>0$.\\
By the canonical bundle formula for blow-ups we have
\begin{equation}\label{canonico}
-K_{\tl X}= -\sigma^\ast K_X - (n-s-1)E = r_X\ac - (n-s-1)E.
\end{equation}
Clearly, $-K_{\tl X}\cdot C_\sigma >0$.\par
\medskip
If $r_X \ge \left[\frac{n}{2}\right]+1$, we get $-K_{\tl X} \cdot \ell =r_X -n +s +1 >0$.
By the Kleiman criterion it follows that $-K_{\tl X}$ is ample, so $\tl X$ is a Fano manifold.
We also get that the length of the ray contracted by $\f$ is $r_X -n +s +1$.\par
\medskip
Assume now that $r_X \ge \left[\frac{n+1}{2}\right]+1$. From (\ref{canonico}) it follows that the line bundle
$$-K_{\tl X} -E = r_X\ac -(n-s)E $$
is ample on $\tl X$, since $(-K_{\tl X} -E)\cdot C_\sigma =n-s$ and
$(-K_{\tl X} -E)\cdot \ell = r_X + s -n>0$.\\
Therefore its restriction to $E$, which by adjunction is $-K_E$, is ample, too; hence $E$ is a Fano manifold.\\
Let $m$ be a positive integer and denote by $D_m$ the divisor $-K_{\tl X} -E +m(\ac -E)$.
Then $D_m$ is ample on $\tl X$, being the sum of an ample line bundle and a nef one,
so $h^1(m\ac-(m+1)E)=h^1(K_{\tl X} +D_m)=0$, by the Kodaira Vanishing theorem.
It follows that the morphism
$$H^0(\tl X, m(\ac-E)) \longrightarrow H^0(E, m(\ac-E)|_E)=H^0(E, m\xi_{N^\ast_{\Sigma/X}(1)})$$
is surjective and  the last claim is proved.
\end{proof}

\begin{remark}\label{gcuse} In the setting of the previous proposition, assume that $r_X \ge [(n+1)/2]+1$, that  $|\xi_{N^\ast_{\Sigma/X}(1)}|$ gives a morphism $\f_E: E \to T$ onto a normal variety and that $\f$ is of fiber type. Then we can assume that $\f$ is defined by the linear system $|\ac-E|$.
\end{remark}

\begin{proof}
Let $\tl \f:\tl X \to \tl Y$ be the morphism defined by the linear system $|\ac-E|$.
Since the fibers of $\tl \f$ are connected, in the Stein factorization of $\tl \f$ the finite morphism $g$ is the normalization (Cf. \cite[1.13]{De}):
$$
\xymatrix{\tl X \ar[rr]^{\tl \f} \ar[rd]_{\f}& & \tl Y \\ & Y \ar[ru]_{g} & }
$$
The divisor $E$ is $\tl \f$-ample, so the restriction of $\tl \f$ to $E$ is onto $\tl Y$; since this restriction, by the last claim of Proposition (\ref{gc}) is $\f_E$  we have $\tl Y = T$, hence $\tl Y$ is normal and $g$ is an isomorphism.
\end{proof}

The following lemma shows that we can apply our construction to manifolds of
Picard number one containing a large linear space
whose normal bundle is nume\-ri\-cally effective.

\begin{lemma}\label{gca}
Let $X \subset \proj^N$ be a smooth variety of dimension $n$ and Picard number one,
containing a linear space  $\Lambda$ of dimension $s$.
Assume that  the normal bundle $N_{\Lambda/X}$ of $\Lambda$ in $X$ is nef
and denote by $c$ the non negative integer such that $\det N_{\Lambda/X}=\Ol_\Lambda(c)$.
Then $X$ is a Fano manifold of index $r_X = s+1+c$ covered by lines.
\end{lemma}

\setcounter{equation}{0}

\begin{proof}
By the adjunction formula we have
$$K_\Lambda=(K_X + \det N_{\Lambda/X})|_{\Lambda},$$
whence
$$(-K_X)|_{\Lambda} = \Ol_\Lambda(s+1+c),$$
from which we can derive
\begin{equation}\label{cano}
-K_X = \Ol_X(s+1+c).
\end{equation}
Let $l$ be a general line in $\Lambda$. From the nefness of $N_{\Lambda/X}$ and the exact sequence
$$\shse{N_{l/\Lambda}= \Ol_{\Lambda}(1) ^{\oplus (s-1)}}{N_{l/X} }{(N_{\Lambda/X})|_l},$$
we have that $N_{l/X}$ is nef.\\
Therefore $l$ is a free rational curve in $X$
(see \cite[Definition II.3.1]{Kob}), which is thus covered by lines
by \cite[Proposition II.3.10]{Kob}.
\end{proof}

\section{Projective $n$-folds covered by linear subspaces of dimension $\ge \frac{n}{2}$}

In this section we will prove that, if a smooth complex projective variety
$X \subset \proj^N$ of Picard number one and dimension $2s$
contains a linear subspace $\Lambda \simeq  \proj^s$ whose normal bundle is $T_{\proj^s}(-1)$, then
$X$ is the Grassmannian of lines in $\proj^{s+1}$.\\
This result, as explained in Corollary
(\ref{exsato}), completes \cite[Main Theorem]{SatoTo},
in which  smooth projective varieties of dimension $n$ covered by linear subspaces
of dimension greater than or equal to $\frac{n}{2}$ were classified.

\begin{lemma}\label{emet} Let $X \subset \proj^N$ be a smooth variety
containing a linear subspace $\Lambda \simeq  \proj^s$
whose normal bundle $N_{\Lambda/X}$ is globally generated and such that
 $h^1(N_{\Lambda/X})=0$. Then $X$ is covered by linear subspaces of dimension $s$.
\end{lemma}

\begin{proof}
The Hilbert scheme of $s$-planes in $X$ is smooth at the point $\lambda$ corresponding to
$\Lambda$.
Let $T$ be the unique irreducible component of the Hilbert scheme containing $\lambda$
and let $Z$ be the universal family; we have the following diagram:
$$\xymatrix{& Z \ar[ld]_q  \ar[rd]^p &\\ T & & X.}$$
Let $z$ be a point in $\Lambda':=q^{-1}(\lambda)$; we consider the differential of $p$ at that point
$$d_zp\colon T_zZ \longrightarrow T_{p(z)}X;$$
this map is the identity when restricted to $T_z\Lambda'$.
Recalling that  $T_\lambda T \simeq H^0(N_{\Lambda/X})$ and considering the exact sequence of the
normal bundle of $\Lambda$ in $X$ we get the following commutative diagram
$$\xymatrix{0 \ar[r] & T_z \Lambda' \ar[r] \ar[d]^{Id} & T_zZ \ar[r] \ar[d]^{d_zp} &
H^0(N_{\Lambda/X}) \ar[r] \ar[d]^{ev}& 0 \\
0 \ar[r] & T_{p(z)} \Lambda \ar[r] & T_{p(z)}X \ar[r]& (N_{\Lambda/X})_{p(z)} \ar[r]& 0}$$
which shows that $d_zp$ is surjective - $ev$ is surjective by the spannedness of $N_{\Lambda/X}$ -
hence $p$ is smooth at $z$.
\end{proof}

\begin{proposition}\label{grass} Let $X \subset \proj^N$ be a smooth variety of Picard number one and dimension $2s$
containing a linear subspace $\Lambda \simeq  \proj^s$ whose normal bundle is $T_{\proj^s}(-1)$.
Then $X$ is the Grassmannian of lines $\G(1,s+1)$.
\end{proposition}

\begin{proof} First of all notice that, by Lemma (\ref{gca}), $X$ is a Fano manifold of index
$r_X= s+1+c = s+2$ covered by lines.\\
Let $\sigma \colon \tl X \to X$ be the blow-up of $X$ along $\Lambda$, denote by
$E=\proj(N_{\Lambda/X}^\ast)$ the exceptional divisor and by $\ac$ the pull-back $\sigma^*\Ol_X(1)$.
By Proposition (\ref{gc}), $\tl X$ is a Fano manifold with a contraction $\f \colon \tl X \to Y$
whose restriction to $E$ is the map associated on $E$ to the linear system
$|m\xi_{N^\ast_{\Lambda/X}(1)}|=|m\xi_{\Omega(2)}|$. This map is, up to a Veronese embedding of the target, the $\pu$-bundle over $\G(1,s)$
given by the projectivization of the universal quotient bundle $\qu$ over $\G(1,s)$,
as shown in Example (\ref{ex1}).\\
Moreover, by Proposition (\ref{gc}), $\cone(\tl X)=\twospan{[C_\sigma]}{[\ell]}$, where
$C_\sigma$ is a rational curve of minimal degree contracted by $\sigma$ and $\ell$ is the strict
transform of a line meeting $\Lambda$ at one point and the length of the extremal ray generated by
$[\ell]$ is $r_X-n+s+1=3$.\\
By Proposition (\ref{fiberlocus}) every non-trivial fiber of the contraction $\f$
has dimension at least two.
Since  $E \cdot \ell=1$ we have that $E$ meets every non-trivial fiber of $ \f$.
As $ \f|_{E}$ is equidimensional with one-dimensional fibers,
it follows that $\f$ cannot have fibers of dimension greater than two, otherwise their
intersection with $E$ would be a fiber of dimension at least two of $\f|_E$.\\
Therefore every non-trivial fiber of $\f$ has dimension two and so, by
Proposition (\ref{fiberlocus}), $\f$ is of fiber type. By Remark (\ref{gcuse}) we can assume that $\f$ is defined by the linear system $|\ac -E|$.\\
Let $F$ be a general fiber of $ \f$; by adjunction we have
$$-K_F = (-K_{\tl X})|_F= ((s+2)\ac -(s-1)E)|_F = 3\ac|_F,$$
hence $(F,\ac|_F) \simeq (\pd, \Ol_\pd(1))$.\\
The line bundle $2\ac -E$ is ample and $(2\ac -E)|_F \simeq \Ol_\pd(1)$;
thus we can apply \cite[Lemma 2.12]{Fuj4} to obtain that $\f$ is a projective
bundle over $\G(1,s)$.\\
Let $\E:=\f_\ast \ac$; the inclusion $E = \proj_{\G(1,s)}(\qu) \hookrightarrow
\tl X = \proj_{\G(1,s)}(\E)$ gives an exact sequence of vector bundles
over $\G(1,s)$
$$\shse{L}{\E}{\qu}.$$
We can compute, using the canonical bundle formula for projective bundles, that
$\det \E= \Ol_{\G(1,s)}(2)$, so, recalling that $\det \qu =\Ol_{\G(1,s)}(1)$,
we have $L=\Ol_{\G(1,s)}(1)$.\\ Since $h^1(\qu^\ast (1))=h^1(\qu)=0$, the sequence splits and
$\tl X= \proj_{\G(1,s)}(\qu \oplus \Ol_{\G(1,s)}(1))$.\\
We have thus proved that the existence in $X$
of a linear subspace $\Lambda \simeq  \proj^s$ whose normal bundle is $T_{\proj^s}(-1)$ completely
determines $X$. As the Grassmannian of lines $\G(1,s+1)$ contains such a linear space
- take  a linear space corresponding to the lines passing through a fixed point -
the proposition is proved.
\end{proof}

\begin{corollary}\label{exsato} $($cf. \cite[Main theorem]{SatoTo}$)$ Let $X \subset \proj^N$ be a smooth
complex variety of dimension $n\ge 2$ covered by linear subspaces of dimensions $s \ge \frac{n}{2}$.
Then $X$ is one of the following:
\begin{enumerate}
\item a $\proj^r$-bundle over a smooth variety. ($r \ge s$);
\item a smooth hyperquadric $\mathbb Q^{2s}$;
\item the Grassmannian of lines $\G(1,s+1)$.
\end{enumerate}
\end{corollary}

\begin{proof} In \cite{SatoTo} the author first showed that the normal bundle
of a general linear subspace is one of the following:
\begin{itemize}
\item[(i)] $\Ol_{\proj^s}^{\oplus a} \oplus \Ol_{\proj^s}(1)^{\oplus (n-s-a)}$;
\item[(ii)] $\Omega_{\proj^s}(2)$;
\item[(iii)] $T_{\proj^s}(-1)$.
\end{itemize}
Then he showed that in case (i) $X$ is a $\proj^r$-bundle over a smooth variety ($r \ge s$) and in case
(ii) $X$ is a smooth hyperquadric.\\
As for case (iii) he showed that $X$ is the Grassmannian of lines in $\proj^{s+1}$, if $s$
is even or if one assumes that all the linear subspaces of the covering family
have normal bundle $T_{\proj^s}(-1)$.\\
Thus to prove the statement it is enough to show that, in case (iii), $X$ is
the Grassmannian of lines in $\proj^{s+1}$; this will follow from Proposition
(\ref{grass}) once we prove that, if $X$ is as in case (iii), then its Picard number
is one.\\
Assume that this is not the case. By Theorem (\ref{nc}) there is an elementary contraction which contracts a covering family of linear subspaces of dimension $s$.
A general fiber of $F$ has dimension at most $2s -1$ and  is covered by linear spaces of dimension $s$; applying \cite[Corollary I.2.20]{Zak}
as in \cite[Theorem 2]{Fuebasta} we derive that $F$ is a projective space,
hence the normal bundle of a general $s$-plane cannot be $T_{\proj^s}(-1)$.
\end{proof}

\begin{corollary}\label{nebastauno}
Let $X \subset \proj^N$ be a smooth variety of  dimension $2s$
containing a linear subspace $\Lambda \simeq  \proj^s$ whose normal bundle is $T_{\proj^s}(-1)$.
Then $X$ is the Grassmannian of lines $\G(1,s+1)$.
\end{corollary}

\begin{proof} By Lemma (\ref{emet}) through every point
of $X$ there is a linear subspace of dimension $s$. By Corollary (\ref{exsato})
$X$ is a $\proj^r$-bundle over a smooth variety, a smooth hyperquadric $\mathbb Q^{2s}$ or
the Grassmannian of lines $\G(1,s+1)$. The first two cases are ruled out since the
corresponding manifolds do not contain a linear subspace with normal bundle
$T_{\proj^s}(-1)$.
\end{proof}

We can now prove the part of Theorem (\ref{main}) regarding manifolds with Picard number greater than one:

\begin{corollary}\label{p2}
 Let $X \subset \proj^N$ be a smooth variety of dimension $2s+1$ and Picard
number greater than one,
containing a linear subspace $\Lambda$ of dimension $s$, whose normal bundle $N_{\Lambda/X}$ is nef.
Assume that $\rho_X >1$; then there is an elementary contraction $\f \colon X \to Y$ which contracts $\Lambda$ and one of the following occurs:
\begin{itemize}
\item $\f$ is a scroll;
\item $Y$ is a smooth curve, and the general fiber of $\f$ is
\begin{itemize}
\item the Grassmannian of lines $\G(1,s+1)$;
\item a smooth hyperquadric $\mathbb Q^{2s}$;
\item a product of projective spaces $\proj^s \times \proj^s$.
\end{itemize}
\end{itemize}
\end{corollary}

\begin{proof} Combine Theorem (\ref{nc}) with Corollary (\ref{exsato}).\end{proof}

\begin{remark} If $\f$ is a scroll and $\dim F \ge s+1$, then $X$ has a projective bundle structure
over $Y$ by \cite[Theorem 1.7]{Ein}. By \cite[Conjecture 14.1.10]{BSbook} this should be the
case also if $\dim F = s$.
\end{remark}

\section{Manifolds with Picard number one - Normal bundles}

Let $X \subset \proj^N$ be a smooth variety of dimension
$2s+1$ and Picard number one, containing a linear subspace $\Lambda$ of dimension $s$
whose normal bundle $N_{\Lambda/X}$ is numerically effective.\\
In this section we will start the proof of the first part of Theorem (\ref{main}), giving the list of
all possible normal bundles of the linear subspace $\Lambda$, showing that $X$ is covered by linear subspaces of dimension $s$ and settling the case  of decomposable normal bundles.

\begin{proposition}\label{swept} Let $X \subset \proj^N$ be a smooth variety of dimension
$2s+1$, containing a linear subspace $\Lambda$ of dimension $s$
whose normal bundle $N_{\Lambda/X}$ is nef. Then $N_{\Lambda/X}$ is one of the following:
\begin{enumerate}
\item $\Omega_\Lambda(2) \oplus \Ol_\Lambda$;
\item $\Omega_\Lambda(2) \oplus \Ol_\Lambda(1)$;
\item $T_\Lambda(-1) \oplus \Ol_\Lambda(1)$;
\item $T_\Lambda(-1) \oplus \Ol_\Lambda$;
\item $\Ol_\Lambda(1)^{\oplus c} \oplus \Ol_\Lambda^{\oplus (s+1-c)}$.
\end{enumerate}
Moreover, through every point of $X$ there is a linear subspace of dimension $s$.
\end{proposition}

\begin{proof}
From the exact sequence
$$\shse{N_{\Lambda/X}}{N_{\Lambda/\proj^N} =  \Ol_{\Lambda}(1)^{\oplus (N-s)}}{(N_{X/\proj^N})|_\Lambda}$$
and the nefness of $N_{\Lambda/X}$ we get that the splitting  of $N_{\Lambda/X}$ on lines
in $\Lambda$ is  of type $(0, \dots, 0,1, \dots, 1)$, hence uniform.\\
By the classification of uniform vector bundles of rank $s+1$ on $\proj^s$ given in
\cite{Ell} and \cite{BalUnif}, taking into account the splitting type, we have that
$N_{\Lambda/X}$ is one of the bundles listed in the statement.
Since all these bundles are generated by global sections and have $h^1(N_{\Lambda/X})=0$ the last assertion follows from Lemma (\ref{emet}).
\end{proof}

\begin{question} Let $X \subset \proj^N$ be a smooth variety of dimension
$2s+1$ such that through every point of $X$ there is a linear subspace of dimension $s$.
It is possible to prove, as in \cite{SatoTo}, that the general linear subspace has
a normal bundle which is spanned at the general point.
Is it true that there exists a linear subspace $\Lambda$ with nef normal bundle?
\end{question}

We  recall a general construction (see \cite[Proof of 0.7]{ABW4}):

\begin{lemma}\label{ein} Let $\Lambda \subset X \subset \proj^N$ be a linear space contained in a
smooth projective variety and such that $N_{\Lambda/X} \simeq N' \oplus \Ol_{\Lambda}(1)$ for some
vector bundle $N'$ over $\Lambda$.
Then there exists a smooth hyperplane section $X'$ of $X$ which contains $\Lambda$
and such that $N_{\Lambda/X'}\simeq N'$.
\end{lemma}

\begin{proof} The existence of the smooth hyperplane section follows from \cite[Corollary 1.7.5]{BSbook};
notice that to apply that result, since $N_{\Lambda/X}^*(1) \simeq N'^*(1) \oplus \Ol_{\Lambda}$,
we do not need assumptions on $\dim \Lambda$;
then by the exact sequence
$$\shse{N_{\Lambda/X'}}{N_{\Lambda/X} \simeq N' \oplus \Ol(1)}{\Ol_\Lambda(1)},$$
we obtain the statement on the normal bundle.
\end{proof}

\begin{proposition}\label{trivial}
Let $X \subset \proj^N$ be a smooth variety of dimension $n \ge 4$
containing a linear space $\Lambda$ of dimension $s$ with $\left[\frac{n}{2}\right] \le s \le n-2$.
Assume that the normal bundle $N_{\Lambda/X}$ is trivial. Then the Picard number of $X$ is at least two.
\end{proposition}

\begin{proof}
Assume by contradiction that $\rho_X=1$; by Lemma (\ref{gca}), $X$ is a Fano manifold
of index $r_X=s+1$ covered by lines. By the first part of Proposition (\ref{gc}),
the blow-up of $X$ along $\Lambda$, which we will denote by $\tl X$, is a Fano manifold with $\rho_{\tl X}=2$,
whose ``other" contraction $\f:\tl X \to Y$ is given by the linear system $|m(\ac -E)|$, where
$\ac:=\sigma^\ast \Ol_X(1)$.
The restriction of $m(\ac -E)$ to $E=\Lambda \times \proj^{n-s-1}$ is
$m\xi_{N^\ast_{\Sigma/X}(1)}=\Ol_E(m,1)$. In particular no curves of $E$ are contracted
by $\f$.\\
The extremal ray associated with $ \f$ is generated by the class $[\ell]$ of the strict transform
of a line meeting $\Lambda$ at one point, hence $E \cdot \ell=1$.
Since $E$ has positive intersection number with curves contracted by $\f$,
it follows that every non-trivial fiber of $\f$ has dimension one.\\
By \cite[Theorem 1.2]{Wicon} $\f$ is either a conic bundle, or a  blow-up of a smooth subvariety of codimension two and in both cases 
$Y$ is smooth.
Assume that $\f$ is a conic bundle; the finite morphism $\f_{|E}\colon E \to Y$
is either birational (if $\f$ has no reducible fibers), or of degree two; since $\rho_Y=1$, in both cases 
we should have $\rho_E=1$. This is clear in the first case, while in the second it follows from \cite{Cor}. So we get a contradiction.\par
\smallskip
If  $\f$ is a blow-up of a smooth
codimension two subvariety then 
the divisor $D=(m+s+1)(\ac -E) = \f^\ast \Ol_Y(m+s+2)$ is nef and big on $\tl X$ for $m>0$; moreover the length of $\f$ is one, so, from Proposition (\ref{gc}) we get $r_X=n-s$, hence $n=2s+1$.\\
This implies that $m\ac-(m+1)E= K_{\tl X} +D$, hence
 $h^1(m\ac-(m+1)E)=0$ by the Kawamata-Viehveg vanishing Theorem; it follows that the morphism
$$H^0(\tl X, m(\ac-E)) \longrightarrow H^0(E,m (\ac-E)|_E)=H^0(E, m\xi_{N^\ast_{\Sigma/X}(1)})$$
is surjective, so the restriction to $E$ of $\f$ is the morphism given
by the linear system $|m\xi_{N^\ast_{\Sigma/X}(1)}|=|\Ol_E(m,1)|$.
In particular the image of $E$ via $\f$ is a smooth divisor isomorphic
to $\proj^s \times \proj^s$. Since $\rho_{ Y}=1$ this is impossible, by Lefschetz's
Theorem on hyperplane sections.
\end{proof}

\begin{remark}\label{dim3} If $n=3$, then, by Lemma (\ref{gca}), $X$ is a Fano manifold
of index $2$ covered by lines, hence a del Pezzo manifold.
By the classification in \cite[Theorem 8.11]{Fuj4}, a del Pezzo threefold of Picard number one
with very ample fundamental divisor is a cubic hypersurface in $\proj^4$,
the complete intersection of two hyperquadrics in $\proj^5$ o a linear section
of $\G(1,4) \subset \proj^9$ with three general hyperplanes.
\end{remark}

\setcounter{equation}{0}

Proposition (\ref{trivial}) allows us to prove that, if the normal bundle of $\Lambda$
is decomposable, the Picard number of $X$ is one and $s \ge 2$, then $X$ is a linear space.

\begin{corollary}\label{obvious}
Let $X \subset \proj^N$ be a smooth variety of dimension $n$ and Picard number one
containing a linear space  $\Lambda$ of dimension $s \ge 2$ with $\left[\frac{n}{2}\right] \le s \le n-2$.
Assume that the normal bundle $N_{\Lambda/X}$ is $  \Ol_\Lambda(1)^{\oplus c} \oplus \Ol_\Lambda^{\oplus (n-s-c)}$.
Then $c=n-s$ and $X$ is a linear space.
\end{corollary}

\begin{proof} By Proposition (\ref{trivial}) we can assume that $c>0$,
so, by Lemma (\ref{ein}) we can find a smooth hyperplane section
$X'$ of $X$ containing $\Lambda$. Then, as in \cite[Theorem 2]{Fuebasta}, we apply \cite[Corollary I.2.20]{Zak},
which yields that $X'$ is a linear space, so we conclude that $X$ is a linear space, too.
\end{proof}

\section{Manifolds with Picard number one - Classification}

In this section we will consider projective manifolds
of dimension $2s+1$ and Picard number one containing a linear subspace of dimension $s$
with numerically effective normal bundle, proving the following

\begin{theorem} Let $X \subset \proj^N$ be a smooth variety of dimension $2s+1$ and Picard number one,
containing a linear subspace $\Lambda$ of dimension $s$, such that its normal bundle $N_{\Lambda/X}$ is nef.
Then $X$ is one of the following:
\begin{itemize}
\item a linear space $\proj^{2s+1}$;
\item a smooth hyperquadric $\mathbb Q^{2s+1}$;
\item a cubic threefold in $\proj^4$;
\item a complete intersection of two hyperquadrics in $\proj^5$;
\item the intersection of the Grassmannian of lines $\G(1,4) \subset \proj^9$ with three general hyperplanes;
\item a hyperplane section of the Grassmannian of lines $\G(1,s+2)$ in its Pl\"ucker embedding.
\end{itemize}
\end{theorem}

\begin{proof}
First of all notice that, by Lemma (\ref{gca}), $X$ is a Fano manifold of index $r_X= s+1+c$.
Moreover, all possible nef normal bundles $N_{\Lambda/X}$ are listed in Proposition (\ref{swept}).\\
When $N_{\Lambda/X}\simeq \Omega_\Lambda(2) \oplus \Ol_\Lambda$, $X$ is a del Pezzo
manifold with very ample fundamental divisor; hence, by the classification in \cite[Theorem 8.11]{Fuj4},
of degree greater than or equal to three.\\
Recalling that the Picard number of $X$ is one and that $X$ contains lines, by the same
classification we have that the degree of $X$ is at most five.\\
A del Pezzo manifold of degree three and Picard number one is a cubic hypersurface in
$\proj^{n+1}$; on the other hand, by
the exact sequence of normal bundles
$$\shse{N_{\Lambda/X}}{\Ol_\Lambda(1)^{\oplus (s+2)}}{\Ol_\Lambda(3)},$$
we see that we cannot have $N_{\Lambda/X} \simeq \Omega_\Lambda(2) \oplus \Ol_\Lambda$, unless
$s=1$.\\
Again by \cite[Theorem 8.11]{Fuj4}, a del Pezzo manifold of degree four is the complete
intersection of two quadric hypersurfaces.\\
Let us show that also this case is possible only if $s=1$. We owe this remark and its proof
to Andrea Luigi Tironi.\\
Let $\mathcal Q$ and $\mathcal Q'$ be the hyperquadrics such that $X= \mathcal Q \cap \mathcal Q'$,
and let $\mathcal F$ be the pencil generated by $\mathcal Q$ and $\mathcal Q'$; by
\cite[Proposition 2.1]{Reidth} the general quadric in $\mathcal F$ is smooth, so we
can assume that $\mathcal Q$ is smooth.\\
By \cite[Corollary 1.7.5]{BSbook} there is a smooth hyperplane section $\mathcal Q_H$ of $\mathcal Q$ containing
$\Lambda$; by the exact sequence of normal bundles
$$\shse{N_{\Lambda/\mathcal Q_H}}{N_{\Lambda/\mathcal Q}}{\Ol_\Lambda(1)},$$
recalling that $N_{\Lambda/\mathcal Q_H} \simeq \Omega_\Lambda(2) \oplus \Ol_\Lambda(1)$
we have $N_{\Lambda/\mathcal Q} \simeq \Omega_\Lambda(2) \oplus \Ol_\Lambda(1)^{\oplus 2}$.\\
Therefore the exact sequence
$$\shse{N_{\Lambda/X}}{N_{\Lambda/\mathcal Q}}{(N_{X/\mathcal Q})|_{\Lambda}}$$
becomes
$$\shse{\Omega_\Lambda(2) \oplus \Ol_\Lambda}{\Omega_\Lambda(2) \oplus \Ol_\Lambda(1)^{\oplus 2}}
{\Ol_\Lambda(2)}.$$
A computation of the total Chern class shows that this is possible only if $s=1$.\\
Again by \cite[Theorem 8.11]{Fuj4}, a del Pezzo manifold of degree five is a linear section of $\G(1,4)$.\\
If  $N_{\Lambda/X} \simeq \Omega_\Lambda(2) \oplus \Ol_\Lambda(1)$, then $X$ is a smooth hyperquadric
by the Kobayashi--Ochiai Theorem \cite{KOc}.\\
In case $N_{\Lambda/X} \simeq T_\Lambda(-1) \oplus \Ol_\Lambda(1)$, by Lemma (\ref{ein})
there exists a smooth hyperplane section $X'$
of $X$ containing $\Lambda$; moreover the normal bundle of $\Lambda$ in $X'$
is $T_\Lambda(-1)$, hence, by Corollary (\ref{nebastauno}), $X'$ is the Grassmannian of lines $\mathbb G(1,s+1)$.
But, by \cite[Corollary 1.3 and Proposition 2.1]{Fuadim}, $\mathbb G(1,s+1)$ cannot be a
hyperplane section of another manifold, unless $s=2$; in this case $X'$ is a four-dimensional
hyperquadric, hence $X$ is a five-dimensional hyperquadric. Note that, since
$T_{\proj^2}(-1) \simeq \Omega_{\pd}(2)$, this was already part of the previous case.\\
As to the remaining possibility, if the normal bundle is decomposable, then $X$ is a linear space by Corollary (\ref{obvious}) while the more difficult case  $N_{\Lambda/X} \simeq T_\Lambda(-1) \oplus \Ol_\Lambda$ is settled in the next subsection.
\end{proof}

\subsection{ Normal bundle isomorphic to $T_\Lambda(-1) \oplus \Ol_\Lambda$}

We will start by proving that $\Lambda$ belongs to a special
one-dimensional family of linear subspaces of $X$:

\begin{prop}\label{lookingforsigma} Let $X \subset \proj^N$ be a smooth variety
of Picard number one and dimension $2s+1$
containing a linear subspace $\Lambda \simeq  \proj^s$
whose normal bundle is $T_{\proj^s}(-1)\oplus \Ol_{\proj^s}$.\\
Then there is a subvariety $\Sigma \subset X$ such that
$(\Sigma, (\Ol_X(1))|_\Sigma)\simeq (\pu \times \proj^s,\Ol_{\pu \times \proj^s}(1,1))$
which contains $\Lambda$ as a fiber of the first projection.
Moreover $\Sigma = \langle \Sigma \rangle \cap X$.
\end{prop}

\begin{proof} By Lemma (\ref{gca}), $X$ is a Fano manifold of index $s+2$ covered by lines.\\
Let $\sigma \colon \tl X \to X$ be the blow-up of $X$ along $\Lambda$, and denote by
$E=\proj(N^\ast_{\Lambda/X})$ the exceptional divisor.
By Proposition (\ref{gc}), $\tl X$ is a Fano manifold with a contraction $ \f \colon \tl X \to Y$
whose restriction to $E$ is the map associated on $E$ to the linear system
$|m\xi_{N^\ast_{\Lambda/X}(1)}|=|m\xi_{\Omega(2)\oplus \Ol(1)}|$, {\em i.e.}, up to a Veronese embedding of the target, the blow-up of $\G(1,s+1)$
along a sub-Grassmannian $\G(1,s)$ as shown in Example (\ref{gblow}).
By Proposition (\ref{gc}) we also have that the extremal ray associated with $ \f$ is generated by the class
$[\ell]$ of the strict transform of a line $l \subset X$ meeting $\Lambda$ at one point and
has length $2$. Let $\ac$ be the pull-back $\sigma^\ast \Ol_X(1)$. Let $A \in \pic(Y)$ be an ample line bundle;
then, for some $t$, $K_{\tl X}+2(\ac+t\f^*A)$ is a supporting divisor for $ \f$.\\
Since  $E \cdot \ell=1$, we have that $E$ meets every non-trivial fiber of $ \f$.\\
As $\f|_{E}$ is equidimensional with one-dimensional fibers,
it follows that $\f$ cannot have fibers of dimension greater than two, otherwise their
intersection with $E$ would be a fiber of dimension at least two of $\f|_E$.
Therefore every non-trivial fiber of $\f$ has dimension at most two.\par
\medskip
We claim that $\f$ is of fiber type.
Assume by contradiction that $\f$ is birational. Then it
is equidimensional by Proposition (\ref{fiberlocus}). We can apply \cite[Theorem 4.1]{AWDuke}
to get that $Y$ is smooth and $ \f$
is the blow-up of a smooth codimension-three center~$T$.\\
Since $E$ meets every non-trivial fiber of $\f$
we have $T\simeq\G(1,s)$.\\
So $Y$ contains $\f(E)\simeq \G(1,s+1)$ as an effective divisor,
but, since $\rho_Y=1$, this implies that $\G(1,s+1)$ is ample in $Y$;
it thus follows by \cite[Corollary 1.3 and Proposition 2.1]{Fuadim} that $s=2$.\\
Therefore $Y$ is a projective space or a smooth hyperquadric and $T \simeq \pd$.
Using the two different blow-up structures of $\tl X$ we can write
$$4\ac -2E = -K_{\tl X} =- \f^\ast K_Y -2\Exc(\f).$$
Therefore the index of $Y$ is even, so $Y \simeq \proj^5$; but
the blow-up of $\proj^5$ along $\pd$ has one fiber type contraction,
so also this case cannot happen.\par
\medskip
It follows that $\f$ is of fiber type. By Remark (\ref{gcuse}) we can assume that $\f$ is defined by the linear system $|\ac -E|$.\\
 Since $E$ meets every fiber of $\f$
we have $Y=\G(1,s+1)$. The restriction of $\f$ to $E$
is birational, hence the general fiber of $\f$ has dimension one. As already noticed, any fiber of $\f$
has dimension at most two, hence $\f$ is a special \ba scroll, so $\tl X=\proj_{\G(1,s+1)}(\E)$,
where $\E:=\f_\ast \ac$.\\
Combining the canonical bundle formula for $\tl X$ as a blow-up with the canonical bundle formula
for $\tl X$ as a \ba scroll, we get
$$ -s(\ac -E) = K_{\tl X} + 2\ac = \f^\ast (K_{\G(1,s+1)} + \det \E)= \f^\ast
\Ol_{\G(1,s+1)}(-s- 2 +\deg \det \E),$$
hence $\det \E=\Ol_{\G(1,s+1)}(2)$. \\
%In particular the splitting type of $\E$ on a line
%in $\G(1,r+1)$ is $(0,2)$ or $(1,1)$. Moreover, if the splitting type is $(0,2)$ then
%the minimal section is a curve contracted by $\sigma$.
%In particular, if $l$ is a line outside $\G(1,r)$ and the splitting type is $(0,2)$ then the intersection
%of $\proj_l(\E)$ with $E$ is the minimal section (the minimal section is in $E$ and $E \cdot f=1$).\\
Denote by $\Lambda_0$ the section of $\sigma\colon E \to \Lambda$ which corresponds to the surjection\linebreak
$\Omega_{\proj^s}(1)\oplus \Ol_{\proj^s}\to\Ol_{\proj^s}$; the restriction of $\xi_{N^\ast_{\Lambda/X}(1)}$ to
$\Lambda_0$ is $\Ol_{\Lambda_0}(1)$, hence $\Lambda_0$ is mapped to a linear subspace $\Lambda_1$ of $\G(1,s+1)$.\\
Let $l$ be any line in $\Lambda_1$; the line in $\Lambda_0$ mapped to $l$ is a section
corresponding to a surjection $\E|_l \to \Ol_l(1)$; hence $\E|_l \simeq \Ol_l(1)^{\oplus 2}$.
By \cite[Th\'eor\`eme]{EHS} the restriction of $\E$ to $\Lambda_1$ is decomposable:
$\E|_{\Lambda_1} \simeq \Ol_{\Lambda_1}(1)^{\oplus 2}$. \\
So $(\proj_{\Lambda_1}(\E|_{\Lambda_1}), \ac) \simeq (\pu \times \proj^s, \Ol_{\pu \times \proj^s}(1,1))$, and
$\Sigma:=\sigma(\proj_{\Lambda_1}(\E|_{\Lambda_1}))$ is a subvariety such that
$(\Sigma, \Ol_X(1))\simeq (\pu \times \proj^s,\Ol_{\pu \times \proj^s}(1,1))$
which contains $\Lambda$ as a fiber of the first projection.\\
To prove the last assertion note that $\Sigma$ is the base locus of the linear subsystem of
$|\Ol_X(1) \otimes {\mathcal I}_\Lambda|$ given by the pull-back of the linear system
$|\Ol_{\G(1,s+1)}(1) \otimes {\mathcal I}_{\Lambda_1}|$.
\end{proof}

Now we will determine the normal bundle in $X$ of the subvariety $\Sigma$ constructed
in the previous proposition.

\begin{prop}\label{ilnormaleequelchee}
Let $X \subset \proj^N$ be a smooth variety of Picard number one and dimension $2s+1$
containing a linear subspace $\Lambda \simeq  \proj^s$ whose normal bundle is
$T_{\Lambda}(-1)\oplus \Ol_\Lambda$.
Let $\Sigma \subset X$ be as in Proposition $(\ref{lookingforsigma})$.\\
Then $N_{\Sigma/X} \simeq p_1^\ast\Ol_\pu(1) \otimes p_2^\ast T_{\proj^s}(-1)$, where $p_1$ and $p_2$
denote the projections of $\Sigma \simeq \pu \times \proj^s$ onto the factors.
\end{prop}

\begin{proof}
By Lemma (\ref{gca}), $X$ is a Fano manifold of index $s+2$ covered by lines.
Let $\sigma \colon \tl X \to X$ be the blow-up of $X$ along $\Sigma$, and denote by
$E=\proj_\Sigma(N^\ast_{\Sigma/X})$ the exceptional divisor.
By Proposition (\ref{gc}),  $E$ is a Fano manifold.
By adjunction
$$\det N_{\Sigma/X}= K_{\Sigma}-(K_X)|_\Sigma = \Ol_\Sigma(s,1)$$
Let $p\colon  E \to \proj^s$ be the composition of the bundle projection
with $p_2$; the fiber $F_x$ of $p$ over a point $x \in \proj^s$ is $\proj_{l_x}((N^\ast_{\Sigma/X})|_{l_x})$
where $l_x$ is the fiber of $p_2$ over $x$.\\
By adjunction $F_x$ is a Fano manifold, hence, recalling that $c_1((N^\ast_{\Sigma/X})|_{l_x})=-s$,
we have that $(N^\ast_{\Sigma/X})|_{l_x} \simeq  \Ol_{l_x}(-1)^{\oplus s}$.\\
So $N_{\Sigma/X}\otimes p_1^\ast\Ol_\pu(-1)$ is trivial on the fibers of $p_2$, hence
$N_{\Sigma/X}\otimes p_1^\ast\Ol_\pu(-1)=p_2^\ast (\F)$, with $\F$ a vector bundle on $\proj^s$.
In particular
$$(N_{\Sigma/X})|_\Lambda \simeq (N_{\Sigma/X}\otimes p_1^\ast\Ol_\pu(-1))|_\Lambda \simeq (p_2^\ast (\F))|_\Lambda =\F.$$
From the the exact sequence of normal bundles
%$$\shse{N_{\Lambda/\Sigma}}{N_{\Lambda/X}}{(N_{\Sigma/X})|_{\Lambda}}$$
$$\shse{\Ol_{\Lambda}}{T_{\Lambda}(-1)\oplus \Ol_{\Lambda}}{(N_{\Sigma/X})|_{\Lambda}},$$
it follows that $(N_{\Sigma/X})|_{\Lambda}$ is nef; recalling that $c_1((N_{\Sigma/X})|_{\Lambda})=1$
we have that $(N_{\Sigma/X})|_{\Lambda}$ is uniform so, either $(N_{\Sigma/X})|_{\Lambda}$ is decomposable or
$(N_{\Sigma/X})|_{\Lambda} \simeq T_{\proj^s}(-1)$. The first case is not possible, since the sequence would
split. Therefore $\F=T_{\proj^s}(-1)$ and $N_{\Sigma/X} \simeq p_1^\ast\Ol_\pu(1) \otimes p_2^\ast T_{\proj^s}(-1)$.
\end{proof}

Now we prove that the existence of a subvariety $\Sigma$ as above completely determines
the manifold $X$.

\begin{prop}\label{dorando} Let $X \subset \proj^N$ be a smooth variety of Picard number one and
 dimension $2s+1$ containing a subvariety  $\Sigma$  such that $\Sigma = \langle \Sigma \rangle \cap X$ and
$(\Sigma, (\Ol_X(1))|_\Sigma))\simeq (\pu \times \proj^s,\Ol_{\pu \times \proj^s}(1,1))$,
with  $N_{\Sigma/X} \simeq p_1^\ast\Ol_\pu(1) \otimes p_2^\ast T_{\proj^s}(-1)$.\\
Let $\sigma \colon \tl X \to X$ be the blow-up of $X$ along $\Sigma$. Then $\tl X$ is a divisor in the linear
system $|2\xi - \f^\ast \Ol_{\G(1,s)}(1)|$
in $\proj_{\G(1,s)}(\Ol_{\G(1,s)}(1) \oplus \qu^{\oplus 2})$, where $\xi$ denotes the tautological line bundle
and $\f$ the bundle projection.
\end{prop}

\begin{proof}
By Lemma (\ref{gca}) $X$ is a Fano manifold of index $s+2$ covered by lines.
Denote by $E=\proj_\Sigma(N^\ast_{\Sigma/X})$ the exceptional divisor; by Proposition (\ref{gc}),
$\tl X$ and $E$ are Fano manifolds.
Moreover, the ray associated with the extremal contraction $\f\colon \tl X \to Y$ different from $\sigma$
has length three and its restriction to $E$ is the map associated on $E$ to the linear system
$|m\xi_{N^\ast_{\Sigma/X}(1)}|=|m\xi_{p_2^\ast \Omega_{\proj^s}(2)}|$, which is, up to a Veronese immersion of the target, the map described in Example (\ref{qb}).
Denote by $\ac$ the pull-back $\sigma^\ast \Ol_X(1)$;  we can
take $K_{\tl X}+3\ac$ as a supporting divisor for $\f$.\\
Since  $E \cdot \ell=1$ we have that $E$ meets every non-trivial fiber of $\f$.
As $\f|_{E}$ is equidimensional with two-dimensional fibers,
it follows that $\f$ cannot have fibers of dimension greater than three, otherwise their
intersection with $E$ would be a fiber of dimension at least three of $\f|_E$.
Therefore every non-trivial fiber of $\f$ has dimension at most three.\par
\medskip
We claim that $\f$ is of fiber type.
Assume by contradiction that $\f$ is birational; then by Proposition (\ref{fiberlocus}) it
is equidimensional. We can apply \cite[Theorem 4.1]{AWDuke} to get that $Y$ is smooth and $\f$
is the blow-up of $Y$ along a smooth center.\\
Since $E \cdot \ell=1$ the intersection of $E$ with a fiber of $\f$ is a $\pd$,
contradicting the fact that fibers of $\f|_{E}$ are isomorphic
to $\pu \times \pu$.\par
\medskip
It follows that $\f$ is of fiber type. By Remark (\ref{gcuse}) we can assume that $\f$ is defined by the linear system $|\ac -E|$.
Since $E$ meets every fiber of $\f$
we have $Y=\G(1,s)$. The contraction $\f$ is supported by $K_{\tl X}+3\ac$, it is
elementary and equidimensional with three-dimensional fibers, so it is a quadric bundle.\\
Let $\E:=\f_\ast \ac$; $\tl X$ embeds in $P:=\proj_{\G(1,s)}(\E)$ as a divisor of relative degree $2$.\\
Let $\E':=\f_\ast (\ac|_E)$; notice that, as shown in Example (\ref{qb}), $\E' \simeq \qu^{\oplus 2}$.\\
The vector bundle $\E$ has $\E'$ as a quotient.
Indeed, if $x \in \G(1,s)$ is a point and we denote by $F$ and $f$
the fibers of $\f$ and ${\f}|_{E}$ over $x$, we have
that $\E'_x =H^0(\ac|_{f})$ is a quotient of $\E_x= H^0(\ac|_F)$.\\
It follows that there exists an exact sequence on $\G(1,s)$:
$$\shse{\Ol_{\G(1,s)}(a)}{\E}{\qu \oplus \qu}.$$
Call $P'$ the projectivization of $\E'$; since $\tl X|_{P'}=E$ we have, by Example (\ref{qb}), that
$$\tl X= 2\ac - \f^\ast \Ol_{\G(1,s)}(1).$$
Combining the canonical bundle formula for $P$,
$$K_{P} +5\ac=\f^\ast (K_{\G(1,s)}+ \det \E),$$
with the blow-up formula giving the canonical bundle of $\tl X$
$$K_{\tl X}=-(s+2)\ac +(s-1)E$$
and the adjunction formula
$$K_{\tl X}=(K_{P} + \tl X)|_{\tl X},$$
we obtain
\begin{eqnarray*}
-(s+2)\ac +(s-1)E &=&  -5\ac+\f^\ast (K_{\G(1,s)}+ \det \E)+2\ac - \f^\ast \Ol_{\G(1,s)}(1)\\
 &= &-3\ac + \f^\ast\Ol_{\G(1,s)}(-s-2+\deg \det \E)\\
  &= &-3\ac + (-s-2+\deg \det \E)(\ac - E)\\
  &= &(-s-5+\deg \det \E)\ac + (s+2+\deg \det \E)E.
\end{eqnarray*}
It follows that $\deg \det \E=3$, therefore $a=1$; since $h^1(\qu^\ast(1)^{\oplus 2})=
h^1(\qu^{\oplus 2})=0$, the above sequence splits and
$\E= \Ol_{\G(1,s)}(1) \oplus \qu^{\oplus 2}$.
\end{proof}

\begin{cor}\label{finalmente}
Let $X \subset \proj^N$ be a smooth variety of Picard number one and dimension $2s+1$
containing a linear subspace $\Lambda \simeq  \proj^s$ whose normal bundle is $T_{\proj^s}(-1)\oplus \Ol_{\proj^s}$;
then $X$ is a hyperplane section of the Grassmannian of lines \linebreak $\G(1,s+2)$.
\end{cor}

\begin{proof} By Propositions (\ref{lookingforsigma}), (\ref{ilnormaleequelchee}), (\ref{dorando}),
there is only one manifold which contains a linear subspace as in the statement.
A smooth hyperplane section $\G(1,s+2) \cap H$ of the Grassmannian of lines $\G(1,s+2)$ contains such a linear space
- just take the intersection of $H$ with a linear space corresponding to lines passing through a fixed point -
so the statement follows.
\end{proof}

\noindent
\small{{\bf Acknowledgements}. The results about extremality of families of lines in \cite{BSW}
were brought to our consideration by a stimulating series of lectures given by Paltin Ionescu.
We thank him also for drawing our attention to the paper \cite{SatoTo}.\\
We thank Jaros{\l}aw Wi\'sniewski who gave us a suggestion for the proof of Lemma
 (\ref{emet}) and Andrea Luigi Tironi, for his careful reading of a preliminary version
of the paper and for his remarks which allowed us to remove a non-effective case from the statement
of Theorem (\ref{main})}.

\providecommand{\bysame}{\leavevmode\hbox to3em{\hrulefill}\thinspace}
\providecommand{\MR}{\relax\ifhmode\unskip\space\fi MR }
\providecommand{\MRhref}[2]{%
  \href{http://www.ams.org/mathscinet-getitem?mr=#1}{#2}
}
\providecommand{\href}[2]{#2}

\end{document}